\documentclass[journal,transmag]{IEEEtran}

\usepackage[section]{placeins}
\usepackage{amsmath, amssymb}
\usepackage{caption}
\usepackage{graphicx}
\usepackage{amsmath,amssymb}
\usepackage{float}
\usepackage{caption}
\usepackage{subcaption}
\usepackage{relsize}
\usepackage{wrapfig}

\usepackage{booktabs}
\usepackage[table]{xcolor}
\usepackage{caption}
\usepackage{graphicx}
\usepackage{array}
\usepackage[most]{tcolorbox}
\usepackage{xcolor}

\renewcommand{\d}{\text{\rmfamily d}}

\ifCLASSINFOpdf
\else
\fi
\begin{document}
%
\title{Minimization of eddy currents in permanent magnets of an electric machine with shape derivatives}


\author{\IEEEauthorblockN{Alessio Cesarano\IEEEauthorrefmark{1}, and Peter Gangl\IEEEauthorrefmark{2},}
\IEEEauthorblockA{\IEEEauthorrefmark{1} Institute of Numerical Mathematics - Johannes Kepler University Linz (JKU), 4040 Linz, Austria}
\IEEEauthorblockA{\IEEEauthorrefmark{2}Johann Radon Institute for Computational and Applied Mathematics (RICAM), ÖAW, 4040 Linz, Austria}}

\markboth{Journal of \LaTeX\ Class Files,~Vol.~XXX, No.~XXX, August~2025}%
{Shell \MakeLowercase{\textit{et al.}}: Bare Demo of IEEEtran.cls for IEEE Transactions on Magnetics Journals}
%



\IEEEtitleabstractindextext{%
\begin{abstract}
In this work we deal with the shape optimization of an electric machine considering time-dependent effects such as eddy currents. The considered electric machine is an interior permanent magnet synchronous machine and we minimize the average dissipated power due to the eddy currents in the magnets over a period of time corresponding to a rotation, while at the same time maximizing the average torque. Our approach is based on the computation of the shape derivative which -- beside the computation of a time discretization of time-dependent state problem -- also involves solving a a time discretization of a time-dependent adjoint problem. The challenge of this problem is related to the dependency of each one of the $N$ time steps of the adjoint problem on two different time steps, due to the use of finite difference in the calculation of eddy current losses.
\end{abstract}

}

\maketitle

\IEEEdisplaynontitleabstractindextext

%
\IEEEpeerreviewmaketitle

\section{Introduction}
%
%
%
%

\IEEEPARstart{E}{ddy current} losses contribute to the total losses in an electric machine. With the justification that design strategies that reduce eddy currents are implemented, such as iron lamination so that the eddy currents do not pass from one sheet to the other, see, e.g., \cite{gyselinckDularHenrotte1999}, eddy currents are often neglected in optimal design. Nevertheless, these effects are still present, e.g., when non-laminated designs are chosen \cite{MellakDeuringerMuetze2022} or in permanent magnets, and their accurate computation is of high relevance. Even though the electrical conductivity of the permanent magnets is lower than the electrical conductivity of iron, eddy currents might be more critical in permanent magnets. This is due to the practical challenges in performing some kind of lamination of the magnets and to the fact that permanent magnets are highly sensitive to heat and subject to de-magnetization. Consequently, also modest eddy currents and consequent temperature increase may be problematic, motivating the minimization of the dissipated power due to eddy current losses in the magnets of interior permanent magnet synchronous machines. 
 
Most often, the models used for the simulation of the physics of an electric machine, as well as for its design optimization, do not include the time derivative of the magnetic field and so the effect of eddy currents is neglected. Indeed, it is common to use the magnetostatic approximation of Maxwell's equations for electric machines, due to the relatively low frequencies at play. In the latter case, eddy currents can be calculated in a post-processing step, for example as a selection criterion in the optimization of electric machines with the use of evolutionary algorithms. A greater accuracy in the calculation of eddy current can be achieved with the use of the magneto-quasi-static approximation of Maxwell's equations, which leads to the eddy current problem, described by a PDE of parabolic type.

\begin{figure}[!t]
\centering  
\includegraphics[width=2.5in]{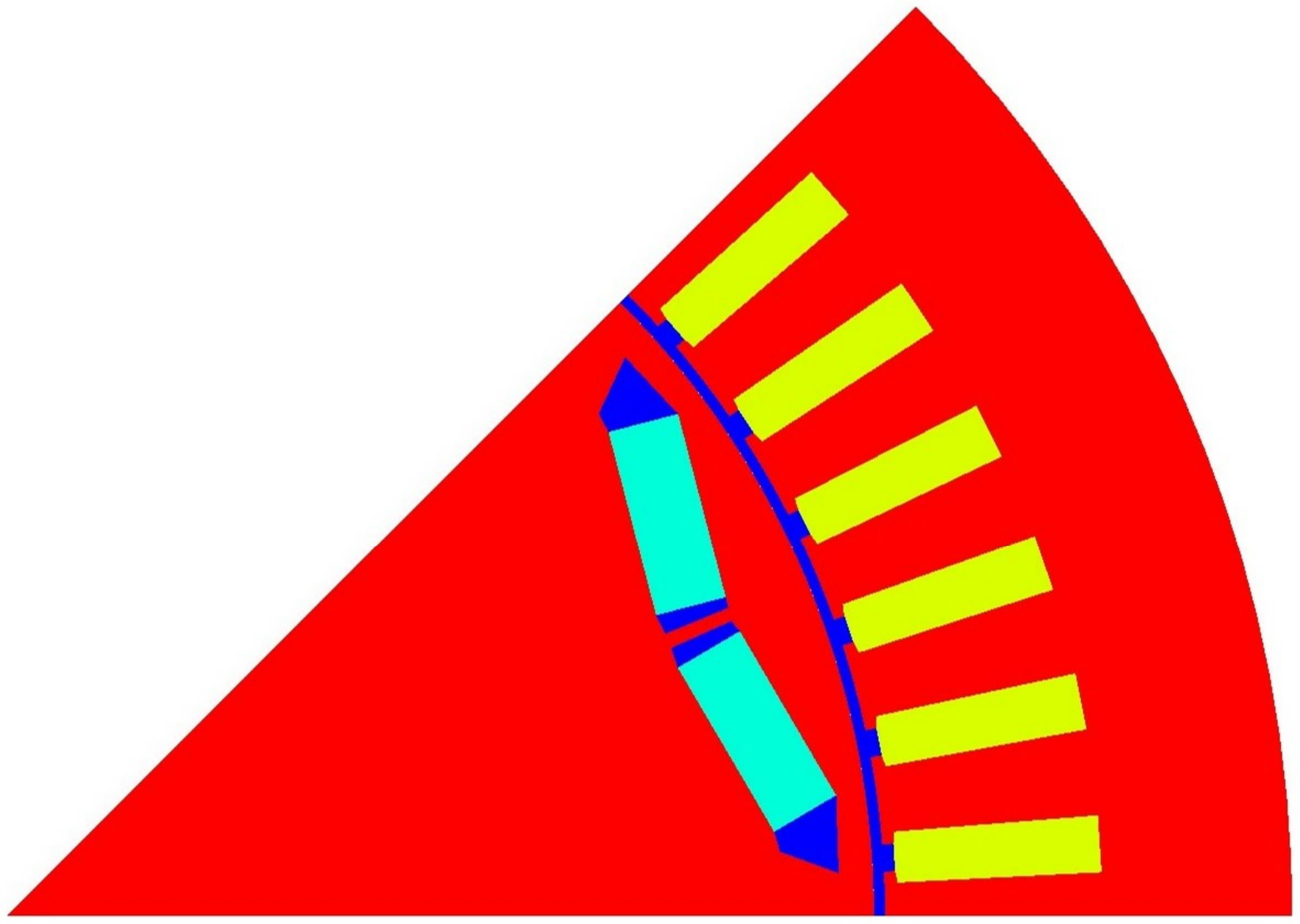}
\caption{Design of the interior permanent magnet synchronous machine (IPMSM) under study. Ferromagnetic material is depicted in red, air is depicted in blue, permanent magnets are depicted in light blue and the coils are depicted in yellow.}
\label{fig.initial_design}
\end{figure}

In this work we solve the eddy current equation in a time-stepping manner and make use of the mathematical notion of shape derivative to implement a gradient-based shape optimization algorithm, with the intent of minimizing eddy currents of an interior permanent magnet synchronous machine (IPMSM), depicted in Fig. \ref{fig.initial_design}. The proposed approach, when compared to evolutionary algorithms, requires much smaller computation times and allows for a larger set of admissible designs.

\section{Eddy currents in permanent magnets} \label{sec_min_dissipated_power}  

\begin{figure}[!ht]
\centering
\subfloat[]{\includegraphics[clip,trim={3.5in} {0.25in} {2.25in} {2.75in},width=1.5in]{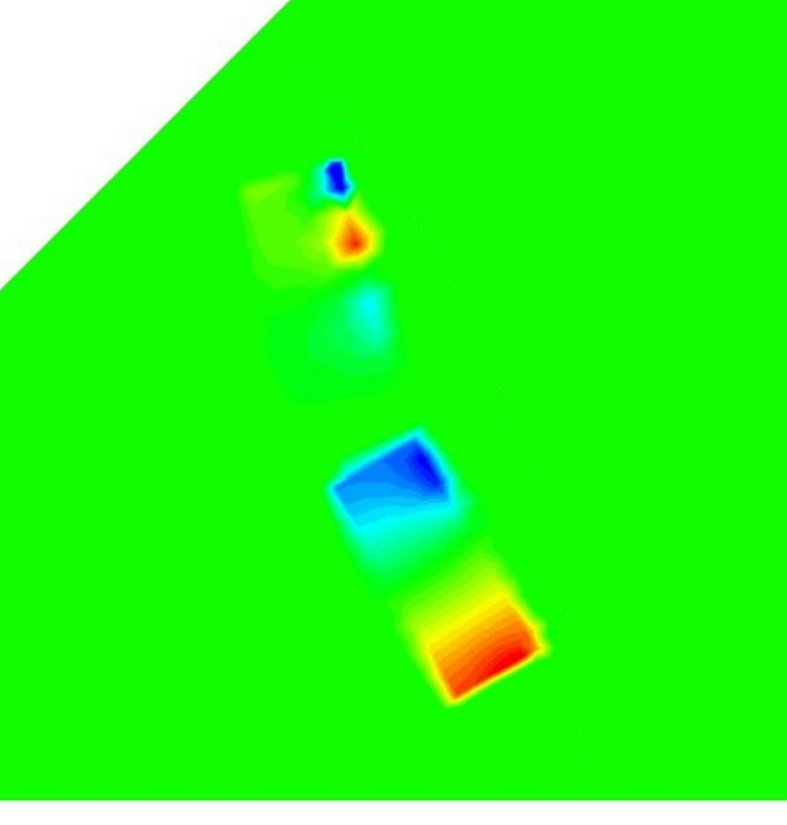}%
\label{fig.jtilde}}
\hfil
\subfloat[]{\includegraphics[clip,trim={3.5in} {0.25in} {2.25in} {2.75in},width=1.5in]{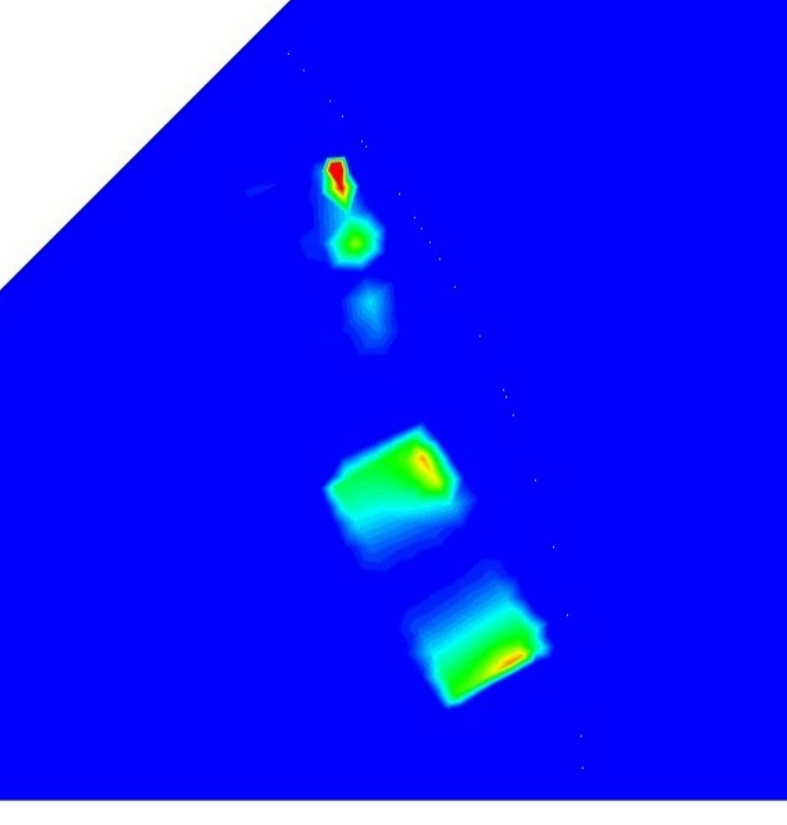}%
\label{fig.p_jtilde}}
\caption{Eddy currents $\tilde{J}_e$ in the permanent magnets at a fixed time step (a). The current is assumed to be purely directed along the z-axis and the loop current to be closed at infinity. Power density related to eddy currents $\tilde{J}_e$ in the permanent magnets at a fixed time step (b).}
\end{figure}

The formula used to describe the eddy currents induced by a variation of magnetic field reads        

\begin{equation} \label{eq_eddycurrents_vec}
    \mathbf{J}_e = - {\sigma} {\partial}_t \mathbf{A},
\end{equation}
where $\mathbf{A}$ is the magnetic vector potential. In the last decades, some work has been devoted to 2-dimensional finite element models that are able to calculate, in an approximate way, the eddy current losses of a corresponding 3-dimensional problem, see, e.g., \cite{SteentjesHameyer2015,Ruoho2009}. In a 3-dimensional domain, eddy currents form a closed loop lying on some plane in the space. Using a 2-dimensional model at a cross section of the machine and introducing a suitable scalar function $u = u(x_1,x_2)$ as the third component of the magnetic vector potential $\mathbf{A}$, one approach to deal with eddy currents is to assume that the eddy currents flow exclusively along the z-axis direction, as depicted in Fig.~\ref{fig.jtilde}. In the latter case, \eqref{eq_eddycurrents} reduces to

\begin{equation} \label{eq_eddycurrents}
J_e = - {\sigma} {\partial}_t u, 
\end{equation} 
 
\noindent where $J_e$ is scalar, but related to a vector directed along the z-axis, meaning that we can directly calculate only axial eddy currents. Thus, the magnets are assumed infinitely long in the axial direction, with the eddy currents loops closed at $z \pm \infty$.  The shorter the axial length with respect to the other dimensions, the worse the approximation. 

As the eddy current paths are closed in one magnet segment, the integral of $J_e$ over the magnet cross section for each rotor position has to be zero. With the subtraction of the average value $\bar{J}_e$ from the calculated current density, as suggested in \cite{Deak2008}, the alternating current density value reads
\begin{equation} \label{eq_eddycurrents_Jtilde}      
    \tilde{J}_e = J_e - \bar{J}_e,
\end{equation}
with
\begin{equation} \label{eq_eddycurrents_Jbar}      
    \bar{J_e} = \frac{1}{|\Omega_e|} \int_{\Omega_e} J_e \; \d x,
\end{equation}
\noindent where $\Omega_e \subset D$ is the region in which we wish to calculate the eddy currents, in our case the permanent magnets region.

%

\section{Magneto-quasi-static problem} \label{sec_seq_2dmagnetostatics}

A necessary step prior to the calculation of the desired cost function and to the shape optimization procedure, is the solution of the state equation. For evaluating the eddy current losses, we are interested in the magneto-quasi-static problem

\begin{align} \label{eq.pbu_sequence_strong}
    \begin{aligned}  
    \sigma \partial_t u 
    - \nabla \cdot \left( \nu(x,|\nabla u|)\nabla u\right) =& f - \nabla \cdot \left(\mathbf{M}^\perp \right), \; \mbox{ in} \: D, \\
    u =& \: 0 \; \qquad \qquad \qquad \mbox{   on }  
    \partial D ,\\
    u =& u_0 \qquad \qquad \qquad \mbox{on } \{0\} \times D,
    \end{aligned}
\end{align}

\noindent where $\nu$ is the magnetic reluctivity, which is nonlinear in the ferromagnetic material, $\sigma$ is the electric conductivity, $\mathbf{M}$ is the magnetization lying in the $xy$-plane and nonzero only in the magnets region, while $f$ is the source term corresponding to the impressed currents in the coils region, that is piecewise constant at a fixed time and evolving in time. 

The PDE in \eqref{eq.pbu_sequence_strong} is also referred to as eddy current equation, which is a mixed parabolic-elliptic PDE, parabolic in the regions where $\sigma$ is nonzero, in our case only on the permanent magnets region, and elliptic elsewhere, see \cite{Cesarano2024}.

Using a finite difference in time, the $j$-th time step can be computed by the variational problem

\begin{align} \label{eq.pbu_sequence_weak}
\begin{aligned}
& \mbox{Find} \:  u_j \in H^1_0(\Omega) \: \text{such that}, \forall w_j \in H^1_0(\Omega), \\
     \int_{\Omega_e} & \sigma_e \: \frac{u_j}{\tau} w_j + \int_{\Omega} \nu \nabla u_j \cdot \nabla w_j = \int_\Omega f_j w_j + \int_{\Omega_e} \sigma_e \: \frac{u_{j-1}}{\tau} w_j,
\end{aligned}
\end{align} 

\noindent for $j=1, \dots, N$. We solve \eqref{eq.pbu_sequence_weak} for $N$ different rotor positions corresponding to one electric period, over a time interval $(0,T)$. The solution of the previous rotor position $u_{j-1}$ must be given for the solution of the problem at the current rotor position $u_j$. As initial condition $u_0$ we choose the solution of the magnetostatic problem. We note that the behavior of the considered machine is periodic for each quarter and anti-periodic for each eighth, allowing for a reduction of the size of the domain. 

\begin{figure}[!ht]
\centering
\subfloat[]{\includegraphics[clip,trim={0.65in} {0.3in} {0.65in} {0.77in},width=1.7in]{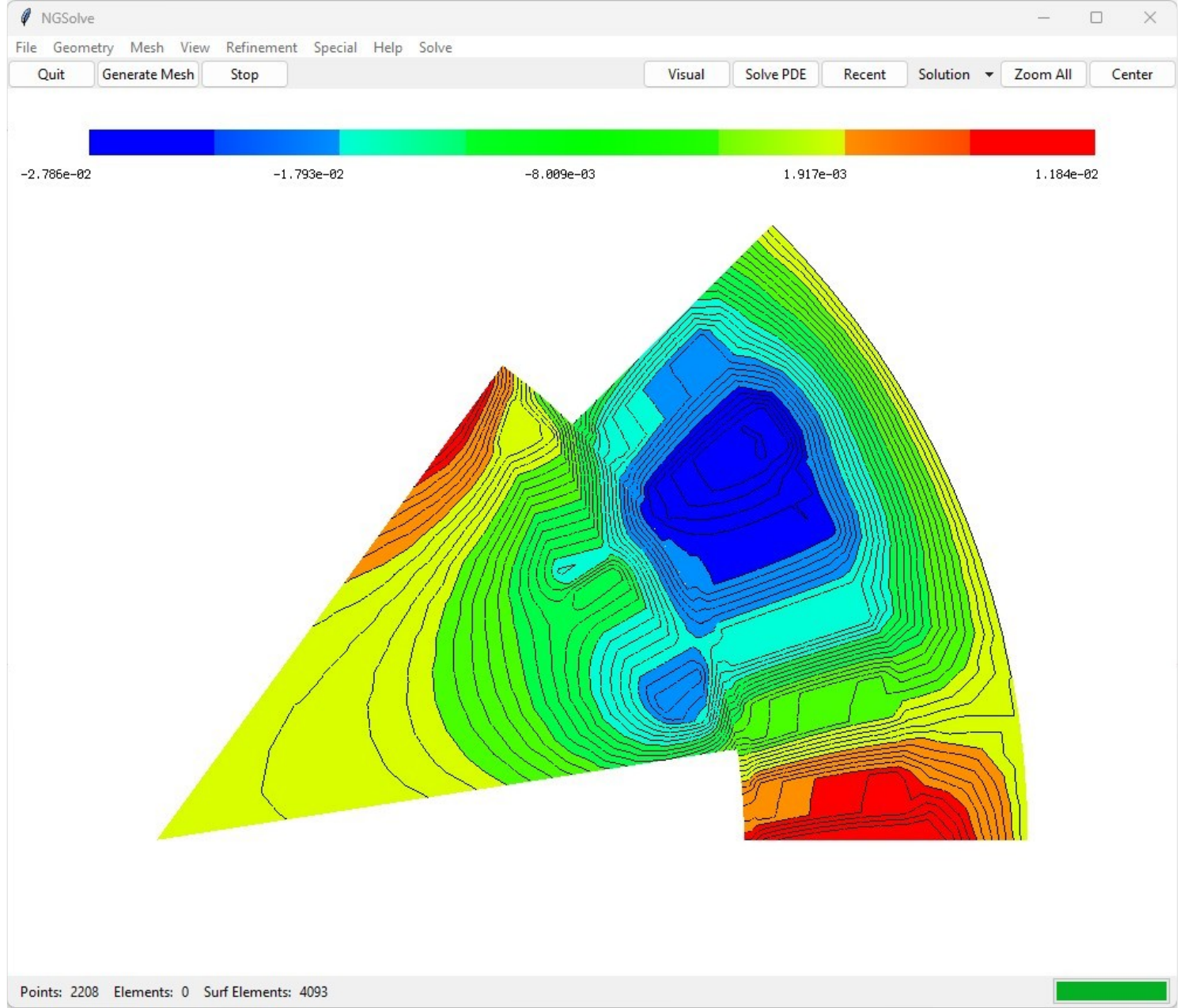}}
\hfil
\subfloat[]{\includegraphics[clip,trim={0.65in} {0.3in} {0.65in} {0.77in},width=1.7in]{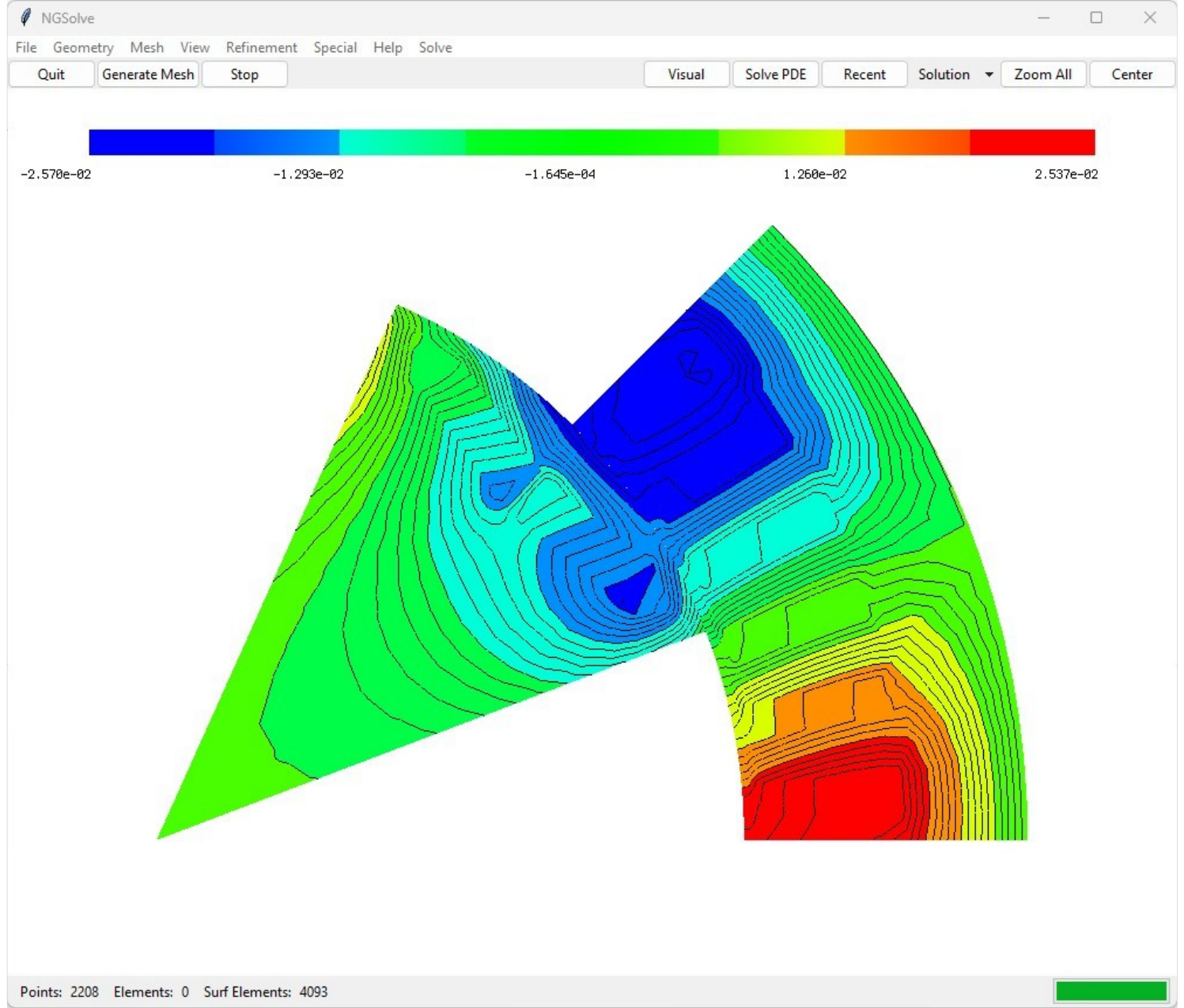}}
\caption{Solution of the problem \eqref{eq.pbu_sequence_weak} at two different step with the locked step method. Antiperiodicity condition is imposed at the left and right sides, as well as at the non connected interface between rotor and stator.}
\label{fig_stateinit}
\end{figure}

Moreover, we use the locked step method for connecting the rotor and the stator of the machine with different angles between each other. The latter consists in identifying algebraically the vertices of the corresponding interfaces of the disjoint rotor and stator domains, after paying attention that the number of vertices is the same in both and that they are equispaced. Different rotor positions are then considered by simply shifting by an integer number such identification. An illustration of the solutions of the problem \eqref{eq.pbu_sequence_weak} at two different rotor positions is depicted in Fig. \ref{fig_stateinit}.


\section{Optimization problem}

In this work, we aim at minimizing the power dissipation due to eddy current losses in the permanent magnets of an IPMSM. In addition, we consider also the average torque, which we intend to maximize, as one of the measures of the performance of the electric motor under study. 

\subsection{Dissipated power due to eddy currents}

The dissipated power related to the eddy current losses according to Joule’s law, at a fixed time step, reads 
\begin{align}\label{eq_dissipated_power}
\begin{aligned}        
    P &= l_z \int_{\Omega_e} \frac{\tilde{J}_e^2}{\sigma} \d x \\
    &= \frac{l_z}{\sigma} \int_{\Omega_e} {\left( J_e - \frac{1}{\Omega_e} \int_{\Omega_e} J_e \; \d y \right)}^2 \d x,
\end{aligned} 
\end{align}

\noindent where we used \eqref{eq_eddycurrents_Jtilde} and \eqref{eq_eddycurrents_Jbar}. The approximation of the time derivative in \eqref{eq_eddycurrents} by finite differences, together with \eqref{eq_dissipated_power}, yields
\small
\begin{align}\label{eq_dissipated_power_2}    
    P_j = \frac{l_z}{\sigma} \int_{\Omega_e} {\left( -  {\sigma} \frac{u_j - u_{j-1}}{\tau}  - \frac{1}{|\Omega_e|} \int_{\Omega_e} - {\sigma} \frac{u_j - u_{j-1}}{\tau} \; \d y \right)}^2 \d x, 
\end{align}
\normalsize 
\noindent where $P_j$ denotes the power dissipation at rotor position $j$, depicted in Fig.~\ref{fig.p_jtilde} for a certain value of $j$, and $\tau$ corresponds to the difference in time between two subsequent rotor positions. We are interested in the average dissipated power, i.e., 
\begin{align} 
P(\underbar{u}) = \frac{1}{N}\sum_{j=1}^{N} P_j.
\end{align}

\subsection{Torque}

For the torque calculation we use a variation of the Maxwell stress tensor method referred to as Arkkio's method, see, e.g., \cite{sadowski1992}, which reads
\begin{align}
T_j = \frac{{\nu}_0 L}{ (r_s - r_r)} \int\limits_\Sigma \, \; ({\nabla u_j})^T \; Q(x,y) \; \nabla u_j \; \mathrm{d} \Sigma, 
\end{align}
where $\Sigma$ is an annulus with inner and outer radii $r_r < r_s$ located in the air gap between the rotor and the stator, $L$ is the length of the machine in the $z$-axis direction and where
\begin{align}
Q(x,y) = \frac{1}{\sqrt{x^2 + y^2}}
\left( \begin{array}{cc} x y & \frac{y^2 - x^2}{2} \\[0.5em] \frac{y^2 - x^2}{2} & -x y \end{array} \right).
\end{align}
As for the dissipated power, we consider the average torque, i.e.,  
\begin{align}
T (\underbar{u}) \; = \; \sum_{j=1}^{N} \; T_j.
\end{align} 

\subsection{Weighted sum of cost functions}

A simple approach to deal with a bi-objective optimization problem consists in converting it into a single-objective optimization problem, performing a scalarization procedure. In this direction, we minimize a convex combination of the given cost functions with fixed, convex parameters, i.e.,
\begin{align} \label{eq_costfunction}
    J(\underbar{u}) = \lambda_1 P(\underbar u) + \lambda_2 (-1) T(\underbar{u}), 
\end{align}
%

\section{Adjoint equations}

We wish to compute the shape derivative of the cost function  \eqref{eq_costfunction}, with the time discretization of the time-dependent PDE \eqref{eq.pbu_sequence_weak} as constraint. In order to find a relevant adjoint equation, let us define the Lagrangian

\begin{align} \label{eq_lagrangian}
\begin{aligned}
    & \mathcal{L}(\Omega, \underbar{u}, \underbar{v}) = \lambda_1 P(\underbar{u}) - \lambda_2 T (\underbar{u}) \: + \\ & \sum_{j=1}^{N} \int_{\Omega_e} \sigma \: \frac{u_j - u_{j-1}}{\tau} v_j \: \d x +\sum_{j=1}^{N} \langle A_j(u_j), v_j \rangle - \langle F_j, v_j \rangle,
\end{aligned}
\end{align}
where $A_j$ are the nonlinear operators corresponding to the magnetostatic problem at the $j$-th rotor position, while $F_j$, $u_j$ and $v_j$ are respectively the corresponding right hand side, solution and test function. 

As $u_j$ is solution of the $j$-th time step of the eddy current problem, the sums are equal to zero and the Lagrangian is equal to the weighted sum of the average dissipated power and average torque, which we want to minimize. 
The time discretization of the adjoint equation is obtained by differentiating the Lagrangian with respect to $u_i$ in the direction of test functions $w_i$, i.e.,

\begin{align}
    0 = \frac{\partial \mathcal{L}}{\partial u_i} (\Omega, \underbar{u}, \underbar{v}) (w_i) = \lambda_1 \frac{\d P}{\d u_i}(\underbar{u})(w_i) - \lambda_2 \frac{\d T}{\d u_i}(\underbar{u})(w_i)\\ + \langle A_i'(u_i)(w_i), v_i \rangle + \int_{\Omega_e}  \: \frac{\sigma}{\tau} ( v_i - v_{i+1} ) \: w_i \: \d x,
\end{align}   
in which we do have some contribution from the sum of the terms coming from the PDEs, and that yields, for each $i \: \in [1,N]$,

\begin{align} \label{eq.pbv_sequence_weak}
\begin{aligned}
& \mbox{Find} \:  v_i \in H^1_0(\Omega) \: \text{such that}, \forall w_i \in H^1_0(\Omega), \\
 & \langle A_i'(u_i)(w_i), v_i \rangle + \int_{\Omega_e}  \: \frac{\sigma}{\tau} \: v_i \: w_i \: \d x   \\ &= - \lambda_1 \frac{\d P}{\d u_i}(\underbar{u})(w_i) + \lambda_2 \frac{\d T}{\d u_i}(\underbar{u})(w_i) + \int_{\Omega_e}  \: \frac{\sigma}{\tau} \: v_{i+1} \: w_i \: \d x.
\end{aligned} 
\end{align} 

The right hand side in \eqref{eq.pbv_sequence_weak} includes the solution of the adjoint in the next rotor position, which must be given. This means that the adjoint equations must be resolved backwards in time, which is expected.       
We note, furthermore, that, as we use order one backward Euler finite differences in time for the calculation of the dissipated power, its derivative with respect to $u_i$ have only 2 non-zero terms out of the $N$ terms related to the different time steps, i.e., 
\begin{equation}\label{eq_adjoint}
    \frac{\d P}{\d u_i}(\underbar{u})(w_i) = \frac1N \frac{\d P_i(u_{i-1},u_i)}{\d u_i}(w_i) + \frac1N \frac{\d P_{i+1}(u_i,u_{i+1})}{\d u_i}(w_i).
\end{equation}  

Taking into account the formula for the dissipated power in \eqref{eq_dissipated_power_2} at the $i$-th and $i+1$-th steps, i.e.,

\begin{align}
\begin{aligned} 
    P_i(u_{i-1}, & u_i) = \\  \frac{l_z}{\sigma} \int_{{\Omega}_e} & {\left( - \sigma \frac{u_i-u_{i-1}}{\tau} - \frac{1}{|{\Omega}_e|} \int_{{\Omega}_e} - \sigma \frac{u_i-u_{i-1}}{\tau} \d y \right)}^2 \d x
\end{aligned}
\end{align}

\noindent and 

\begin{align}
\begin{aligned} 
    P_{i+1}(u_i, & u_{i+1}) = \\ \frac{l_z}{\sigma} \int_{{\Omega}_e} & {\left( - \sigma \frac{u_{i+1}-u_i}{\tau} - \frac{1}{|{\Omega}_e|} \int_{{\Omega}_e} - \sigma \frac{u_{i+1}-u_i}{\tau} \d y \right)}^2 \d x,
\end{aligned}
\end{align} 

\noindent the two terms on the right hand side in \eqref{eq_adjoint} read

\begin{align} \label{eq_rhs_adjoint_bw}
\begin{aligned} 
    \frac{\d P_i(u_{i-1},u_i)}{\d u_i}&(w_i) = \\ 
    \frac{2 l_z}{\sigma} \int_{{\Omega}_e} & \left( - \sigma \frac{u_i-u_{i-1}}{\tau} - \frac{1}{|{\Omega}_e|} \int_{{\Omega}_e} - \sigma \frac{u_i-u_{i-1}}{\tau} \d y \right)  \\ 
    & \cdot \left( - \sigma \frac{w_i}{\tau} - \frac{1}{|{\Omega}_e|} \int_{{\Omega}_e} - \sigma \frac{w_i}{\tau} \d y \right) \d x
\end{aligned}
\end{align}      

\noindent and 

\begin{align} \label{eq_rhs_adjoint_fw}
\begin{aligned} 
   \frac{\d P_{i+1}(u_i,u_{i+1})}{\d u_i}&(w_i) = \\ 
   \frac{2 l_z}{\sigma} \int_{{\Omega}_e} & \left( - \sigma \frac{u_{i+1}-u_i}{\tau} - \frac{1}{|{\Omega}_e|} \int_{{\Omega}_e} - \sigma \frac{u_{i+1}-u_i}{\tau} \d y \right)  \\ 
   & \cdot \left( \sigma \frac{w_i}{\tau} - \frac{1}{|{\Omega}_e|} \int_{{\Omega}_e} \sigma \frac{w_i}{\tau} \d y \right) \d x.
\end{aligned} 
\end{align} 

We remind the reader that $w_i$ is a test function in \eqref{eq_rhs_adjoint_bw} and \eqref{eq_rhs_adjoint_fw}, in which its average is also present. In our numerical implementation, we include this term by modifying a posteriori the assembled finite element matrix arising from \eqref{eq.pbv_sequence_weak}.

\section{Shape optimization framework}
We aim at optimizing the shape of the air pockets adjacent to the permanent magnets in the rotor in order to minimize the aggregated cost function \eqref{eq_costfunction}. We employ the mathematical concept of shape derivatives which, after discretization, can be associated with sensitivities with respect to mesh node positions.

\subsection{Shape derivative}
The shape derivative of a domain-dependent cost function measures its sensitivity with respect to a smooth variation of the domain represented by a vector field $\theta \in H^1(D, \mathbb R^2)$. Defining the transformation $T_t^\theta: D \rightarrow D, x \mapsto = x+t \theta(x)$, the shape derivative of a shape function $\mathcal J = \mathcal J(\Omega)$ is defined by
\begin{align} \label{eq_defdJ}
    d \mathcal J(\Omega; \theta) = \underset{t \rightarrow 0}{\mbox{lim }} \frac{\mathcal J(T_t^\theta(\Omega)) - \mathcal J(\Omega)}{t}
\end{align}
if this limit exists and the mapping $\theta \mapsto d \mathcal J(\Omega; \theta)$ is linear and bounded. Following the approach presented in \cite{Sturm2015}, we can obtain the shape derivative of a PDE-constrained shape optimization problem as
\begin{align} \label{eq_dJ_dtG}
    d \mathcal J(\Omega; \theta) = \partial_t G(0, \underbar{u}, \underbar{v})
\end{align}
with the transformed Lagrangian $G(t, \underbar{u}, \underbar{v}) = \mathcal L(T_t^\theta(\Omega), \underbar{u}\circ (T_t^\theta)^{-1}, \underbar{v}\circ (T_t^\theta)^{-1})$ with $\mathcal L$ as defined in \eqref{eq_lagrangian}, see also \cite{ganglmerkelschoeps2021} for a related study. The differentiation of the right hand side of \eqref{eq_dJ_dtG} can be tedious and error-prone and has been automated within the finite element software package NGSolve in \cite{gangl2021}. In this work, we use the automated shape differentiation capability of NGSolve. For that purpose it is enough to provide the solutions to state and adjoint equations $\underbar{u}$ and $\underbar{v}$, respectively, as well as the Lagrangian function \eqref{eq_lagrangian} in order to obtain the numerical representation of the shape derivative \eqref{eq_defdJ} of problem \eqref{eq_costfunction}.

\subsection{Optimization algorithm}
Given availability of the shape derivative for a given design $\Omega$, the optimization iteratively proceeds as follows: 
\begin{enumerate}
    \item We extract a descent direction, which is a vector field for which the shape derivative is negative, thus yielding a decrease of the cost function for a small enough deformation. This is accomplished by solving a boundary value problem of the form
\begin{align}
\mbox{Find } \theta \in H: \; b(\theta, W) = - d \mathcal J(\Omega; W) \quad \forall W \in H,
\end{align}
where $H \subset H^1(D, \mathbb R^2)$ is a vector-valued finite element space and $b(\cdot, \cdot)$ is a positive bilinear form defined on $H$.
    \item We advect all mesh nodes a small distance $t$ in the direction given by the descent direction $\theta$.

\end{enumerate}
In our numerical implementation, we choose $H$ as the space of all the functions defined on the ferromagnetic and air regions of the rotor, excluding the magnets and the stator, such that every $W \in H$ is zero at the outer boundary of the rotor and the boundaries of the magnets. This choice allows the deformation vector field $\theta$ to act only at the interfaces between the ferromagnetic and air regions of the rotor. As a bilinear form, we use

\begin{equation}\label{eq.bilinear_form_descent}
b(\theta, W) = \frac{1}{2} \langle \nabla \theta + {\nabla \theta}^T, \nabla W + {\nabla W}^T \rangle, 
\end{equation}
which we enrich with a Cauchy-Riemann type term that helps preserving the mesh quality, see, e.g., \cite{IglesiasSturm2018}. As we deform the domain exclusively advecting the mesh, it is important to note that the mesh quality may decrease in the course of the optimization process. For this reason, we assess the mesh quality at each iteration and use it as an additional stopping criterion. In this regard, it could be beneficial to also remesh, see, e.g., \cite{mmg}, but this possibility is not exploited in this work.

\section{Numerical results}

\begin{figure}[!t]
\includegraphics[clip,trim={0in} {0.25in} {0in} {0.25in},width=3.25in]{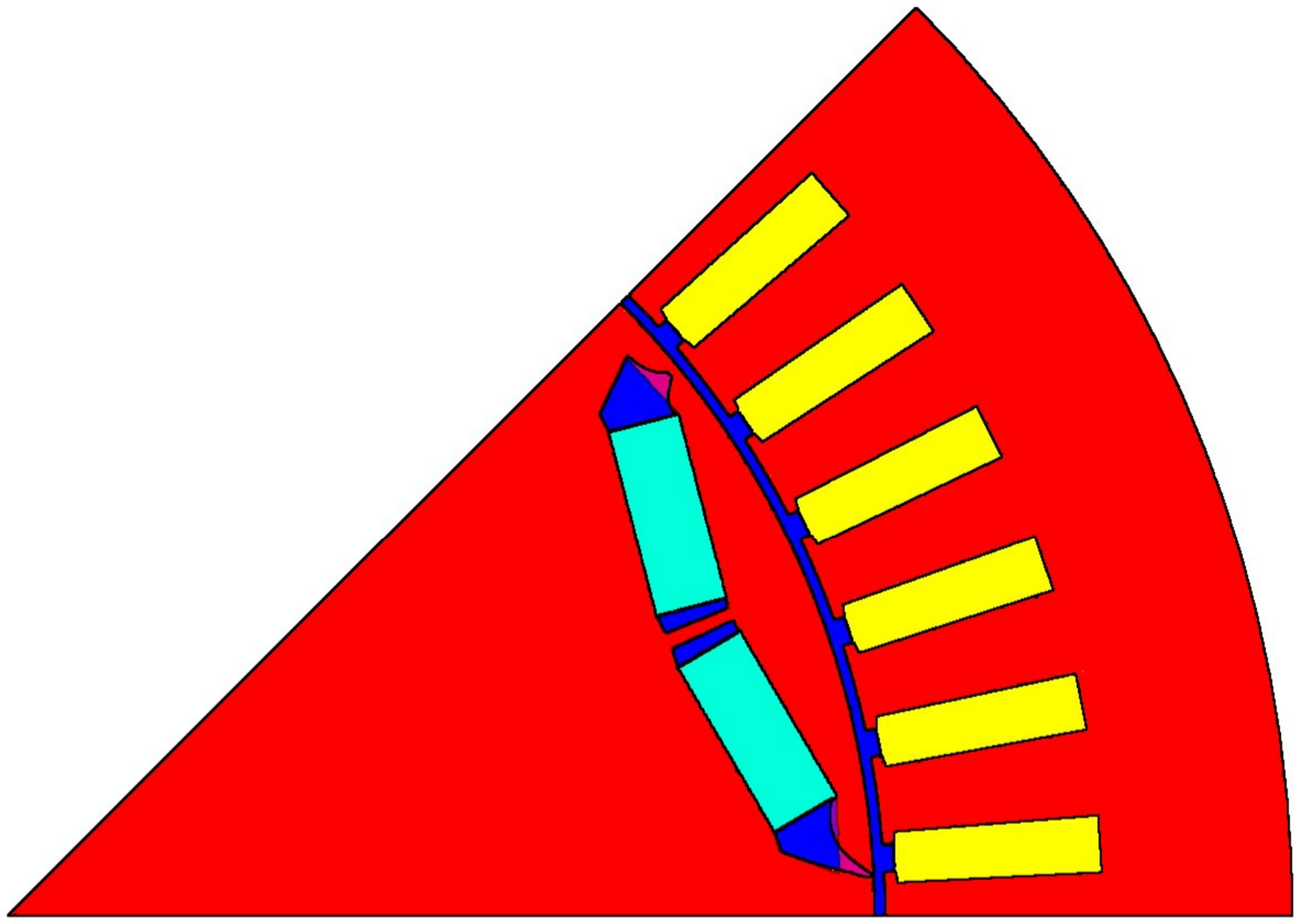}
\caption{Superposition of the initial design of the machine under study and the optimal design obtained with the proposed shape optimization algorithm. }
\label{fig.final_design}
\end{figure}

In this section we present the numerical results of the shape optimization algorithm taking into account the cost functions and the adjoint problem described in the previous sections. We consider the cost function in \eqref{eq_costfunction} with $\lambda_1 = 1 \cdot 10^5$ and $\lambda_2= 1 \cdot 10^{-4}$. The initial design of the machine under scrutiny is depicted in Fig. \ref{fig.initial_design}, while the optimal design obtained with the proposed shape optimization algorithm is depicted in superposition to the initial design in Fig. \ref{fig.final_design}. The average dissipated power related to the eddy currents in the permanent magnets, averaged over an electric period, is reduced by roughly $17\%$, more precisely from $P_{it0}=0.367 \mbox{W}$ to $P_{it53}=0.314 \mbox{W}$, while the average torque increases from $T_{it0}=574.68 \mbox{Nm}$ to $T_{it53}=587.30 \mbox{Nm}$. The optimization process is interrupted due to the approaching change of topology of the design, which is not allowed by our shape optimization algorithm and that leads to the flipping of at least one triangle of the mesh. The solution of the state problem requires about 12 seconds, the solution of the adjoint problem, together with the assembling of the shape derivative requires ca. 3 seconds and the extraction of a descent direction requires, instead, ca. 2 seconds. The results are obtained with $2208$ nodes, an angular velocity of $1500$RPM and $N = 15$ rotor positions.

\begin{figure}[!t]
\includegraphics[clip,trim=0 {1.25in} 0 {1.25in},width=3.5in]{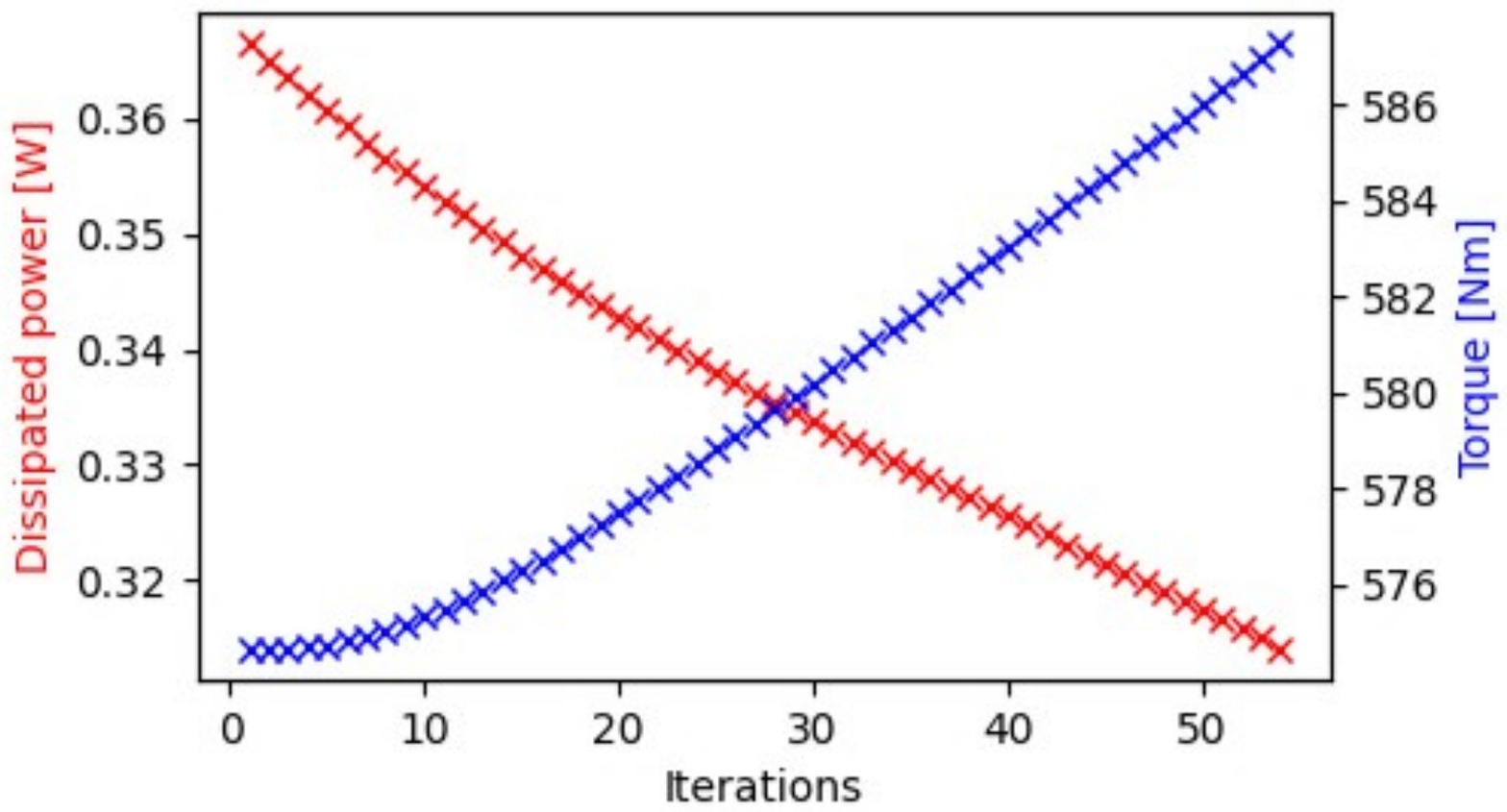}
\caption{Plot of the convergence history of the optimization process. The average dissipated power due to the eddy currents in the magnets over an electric period is reduced from $P_{it0}=0.367 \mbox{W}$ to $P_{it53}=0.314 \mbox{W}$, while the average torque increases from $T_{it0}=574.68 \mbox{Nm}$ to $T_{it53}=587.30 \mbox{Nm}$.}
\label{fig.plot_results}
\end{figure}

\section{Conclusions and future work}

The proposed work deals with the shape optimization of an interior permanent magnet synchronous electric machine, with the consideration of several rotor positions. In particular, the average dissipated power due to eddy currents in the magnets is minimized, while the average torque is maximized at the same time. For a more accurate description of the time-dependent effects of the eddy currents, the time-dependent eddy current equation has been considered. The final design obtained with the proposed shape optimization algorithm presents significant reduction of the average eddy current losses and a greater average torque. However, the latter seems not to be enough mechanically stable to be manufactured without post-processing and the addition of material. For this reason, further developments may include the addition of a mechanical stiffness constraint in the optimization procedure. Future work may also include the solution of the magneto-quasi-static problem with time periodicity.

\section*{Acknowledgment}
The work of A.C. and P.G. is partially supported by the joint DFG/FWF Collaborative Research Centre CREATOR (DFG: Project-ID 492661287/TRR 361; FWF: 10.55776/F90) at TU Darmstadt, TU Graz, JKU Linz and RICAM Linz. P.G. is partially supported by the State of Upper Austria.

\ifCLASSOPTIONcaptionsoff
  \newpage
\fi

\end{document}